\documentclass[11pt,twoside]{article}

\date{}
\usepackage{epsfig}
\begin{document}
\begin{figure}[h]
{\Large{\bf  {\emph{Mathematical Sciences}}}}
\hspace{1.1cm}\hspace{0.9cm}
\hspace{1.1cm}{\small{\bf  {\emph{Vol. 1, No. 1,2 (2007) 01-12}}}}
\end{figure}

\begin{center}
{\Large \bf   Classification of $n-$th order linear ODEs up to
projective transformations}

{\small \bf   Mehdi Nadjafikhah$^{a,}$\footnote{\footnotesize
Corresponding Author. E-mail Address: m$\_$nadjafikhah@iust.ac.ir
}, Seyed Reza Hejazi$^b$}
\end{center}
{\tiny \bf   $^a$Fact. of Math., Dept. of Pure Math., Iran University of Science and Technology, Narmak, Tehran, I.R. Iran.}\\
{\tiny \bf   $^b$Same address.}

\input{amssym}
\newtheorem{Theorem}{\quad Theorem}[section]
\newtheorem{Definition}[Theorem]{\quad Definition}
\newtheorem{Corollary}[Theorem]{\quad Corollary}
\newtheorem{Lemma}[Theorem]{\quad Lemma}
\newtheorem{Example}[Theorem]{\quad Example}
\begin{abstract}
Classification of $n-$th $(n\geq2)$ order linear ODEs is
considered. The equation reduced to \textit{Laguerre–Forsyth} form
by a point transformation then, the other calculations would have
done on this form. This method is due to
\textit{V.A. Yumaguzhin}. 
\\
\copyright {\small \hspace{0.15cm}200x Published by Islamic AZAD
University-Karaj Branch.}
\end{abstract}
\section*{Introduction}
The local classification of linear ODEs up to projective
transformations is obtained in this article. For $n\leq 2$, it is
well known that any $n-$th order linear ODE can be transformed
locally to the form $y^{(n)}=0$ by a point transformation. For
$n\geq0$, this statement is incorrect: there is finite number of
different equivalence classes of linear ODEs.

First this problem was posed by classics of the 15 century E.
Laguerre, G.H. Halphen and others. They obtain results concerning
classification of third and fourth orders linear ODE. Here, this
problem is solved for $n\geq0$ in a neighborhood of regular germs.

Consider a general $n-$th order ODE which is solved by the higher
order derivative
\begin{eqnarray}\label{eq:1}
y^{(n)}=\sum_{i=1}^na_{n-i}(x)y^{(n-i)},
\end{eqnarray}
where $y(x)$ is a smooth function of $x$.

Lie shows that the point symmetry group of a second ordinary
linear differential equation has dimension at most eight,
conversely the equation admits an eight-dimensional symmetry group
if and only if it can be mapped, by a point transformation, to the
linear equation $y''=0$. Thus, the main result is any linear
second ordinary differential equation can mapped to the equation
$y''=0$. So, the condition of second ordinary linear differential
equation is specified.

A same result shows that for $n\geq3$, any linear ODE admits at
most an $(n+4)-$dimensional symmetry group of point
transformation, therefore, the symmetry group is
$(n+4)-$dimensional if and only if the equation is equivalent to
the linear equation $y^{(n)}=0$. In continuation we will work on
the general form of linear ODE in the form of (\ref{eq:1}) once
$n\geq3$.
\section{Laguerre-Forsyth form}
The classification of linear differential equations is a special
case of the general problem of classifying differential operators,
which has a variety of important applications. Consider an $n-$th
order ordinary differential operator corresponding to (\ref{eq:1})
\begin{eqnarray}\label{eq:2}
{\cal D}=a_n(x)D_x^n+a_{n-1}(x)D_x^{n-1}+\cdots+a_1(x)D_x+a_0.
\end{eqnarray}
The aim is finding out when two operators, or two linear ODE, of
type (\ref{eq:2}), can be mapped to each other by a suitable
change of variables. To preserve linearity, we restrict to those
of the form
\begin{eqnarray}\label{eq:3}
\bar{x}=\varphi(x),\qquad\bar{y}=\psi(x)y,
\end{eqnarray}
the chain rule action shows that
$D_{\bar{x}}=(\varphi(x))^{-1}D_x$, and with a rescaling of the
dependent variable by $\displaystyle\psi(x)=e^{\varphi(x)}$ we
obtain the \textit{gauge factor}. So, two differential operator
$\overline{\cal D}$ and $\cal D$ is called \textit{gauge
equivalent} if they satisfy
\begin{eqnarray}\label{eq:4}
{\overline{\cal D}}=\psi\cdot{\cal
D}\cdot\frac{1}{\psi},\qquad\bar{x}=\varphi(x).
\end{eqnarray}
A straightforward calculation shows that the change of variables
(\ref{eq:3}) is given by
\begin{eqnarray*}
\bar{x}=\varphi(x)=\int\frac{dx}{\sqrt[n]{|a_n(x)|}},\;\;\;\;\psi(x)=|a_n(x)|^{\frac{1-n}{2n}}\exp\Big\{\int^x\frac{a_{n-1}(y)}{na_n(y)}dy\Big\},
\end{eqnarray*}
thus (\ref{eq:1}) is gauge equivalent to an operator of the form
\begin{eqnarray}\label{eq:6}
{\cal D}=\pm D_x^n+a_{n-2}(x)D_x^{n-2}+\cdots+a_0(x).
\end{eqnarray}

If $\rho(x)$ be a nonvanishing smooth function, two differential
operator $\overline{\cal D}$ and $\cal D$ is called
\textit{projective equivalence} if they satisfy
\begin{eqnarray}\label{eq:5}
{\overline{\cal D}}=\rho\cdot\psi\cdot{\cal
D}\cdot\frac{1}{\psi},\qquad\bar{x}=\varphi(x).
\end{eqnarray}
A nonsingular $n-$th order linear operator of type (\ref{eq:6}) is
projectively equivalent to one in Laguerre-Forsyth form
\begin{eqnarray}\label{eq:7}
{\cal D}=D_x^n+a_{n-3}(x)D_x^{n-3}+\cdots+a_0(x),
\end{eqnarray}
with change of variable (\ref{eq:5}) in the form of
\begin{eqnarray*}
\bar{x}=\varphi(x),\qquad\bar{y}=\varphi_x^{\frac{n-1}{2}}y,\qquad\rho=\varphi_x^{-n},
\end{eqnarray*}
where $\varphi(x)$ is a solution of the \textit{Schwarzian}
equation
\begin{eqnarray*}
\frac{n(n^2-1)}{12}\frac{\varphi_x\varphi_{xxx}-\frac{3}{2}\varphi_{xx}^2}{\varphi_x^2}=a_{n-2}(x).
\end{eqnarray*}
\section{Classification of linear ODEs of Laguerre-Forsyth form}
A useful theorem help us to reduce the classification of ODEs up
to a special transformation.
\begin{Theorem}
Let $\Delta_1$ and $\Delta_2$ be ODEs of the form (\ref{eq:7}). If
there is a point transformation that takes $\Delta_1$ to
$\Delta_2$, that is
\begin{eqnarray}\label{eq:8}
f(x)=\frac{ax+b}{cx+d},\quad\hat{f}(x,y)=|f'|^{\frac{n-1}{2}}\cdot
y,\quad a,b,c,d\in\bf{R}.
\end{eqnarray}
\end{Theorem}
A transformation $(f,\hat{f})$ of the form (\ref{eq:8}) is
generated by a projective transformation $f$ on $\bf R$. The
isomorphisms $f\rightarrow(f,\hat{f})$ makes a group of point
transformations in the form of (\ref{eq:8}). Consider these
projective transformations in a group $G$ and denoted by all
projective transformations of $\bf R$,i.e.,
\begin{eqnarray*}
G=\Big\{f(x)=\frac{ax+b}{cx+d}\Big|a,b,c,d\in{\bf
R}\;\mbox{and}\;ad\neq bc\Big\}.
\end{eqnarray*}
It is easy to check that $G$ has two connected component
$G_1=\{f\in G|f'>0\}$ and $G_2=\{f\in G|f'<0\}$, thus, $G=G_1\cup
G_2$.
\subsection{Bundles of Laguerre-Forsyth form}
Consider $x$ as a coordinate on $\bf R$ and
$a_{n-3},a_{n-2},...,a_0$ coordinates on ${\bf R}^{n-2}$. Then, we
can construct a fiber bundle corresponding to (\ref{eq:7}) in the
form of
\begin{eqnarray}\label{eq:9}
p:{\bf R}\times{\bf R}^{n-2}\rightarrow{\bf R}.
\end{eqnarray}
Any ODE of type (\ref{eq:7}) identifies with
$\Delta=\{p_n=a_{n-3}(x)p_{n-3}+\cdots+a_0(x)p_0\}$ is a section
of (\ref{eq:9}) denoted by
$S_\Delta:x\rightarrow(x,a_{n-3}(x),...,a_0(x))$, where the
identification $\Delta\rightarrow S_\Delta$ is a bijection.

Let
$\Delta_2=\{\tilde{p}_n=\tilde{a}_{n-3}(\tilde{x})\tilde{p}_{n-3}+\cdots+\tilde{a}_0(\tilde{x})\tilde{p}_0\}$
be an ODE of the form (\ref{eq:7}). Subjecting $\Delta_2$ to an
transformation (\ref{eq:8}), the, we obtain linear ODE
$\Delta_1=\{p_n=a_{n-3}(x)p_{n-3}+\cdots+a_0(x)p_0\}$. The
coefficients $\Delta_2$ are expressed in terms of coefficients of
$\Delta_1$ and projective transformation $f^{-1}$ by the equation
\begin{eqnarray}\label{eq:10}
\tilde{a}_{n-i}=F_{n-i}\Bigg(a_{n-3},...,a_{n-i};\frac{df^{-1}}{d\tilde{x}},...,\frac{d^{i+1}f^{-1}}{d\tilde{x}^{i+1}}\Bigg),\quad
i=3,4,...,n.
\end{eqnarray}
The equation (\ref{eq:10}) is a lifting of a projective
transformation $f$ to diffeomorphism $\bar{f}:{\bf R}\times{\bf
R}^{n-2}\rightarrow{\bf R}\times{\bf R}^{n-2}$ such that
$p\circ\bar{f}=f\circ p$.

For any $f\in G$, a transformation of sections of $p$ defined by
the formula
\begin{eqnarray*}
S\rightarrow f(S)=\bar{f}\circ S\circ f^{-1},
\end{eqnarray*}
then, equation (\ref{eq:10}) can be represented as
$S_{\Delta_2}=f(S_{\Delta_1})$.
\begin{Lemma}\label{lem:1}
Consider two equation of the form (\ref{eq:7}). Then a
transformation $(f,\hat{f})$ of the form (\ref{eq:8}) maps
$\Delta_1$ to $\Delta_2$ if and only if
$S_{\Delta_2}=f(S_{\Delta_1})$.
\end{Lemma}
The main result of the lemma (\ref{lem:1}) is the classification
of ODEs of the form (\ref{eq:7}) up to transformation (\ref{eq:8})
reduces to classification of germs of sections of $p$ up to
projective transformation on $\bf R$.
\subsection{Classification of regular germs}
Let $S$ be a section of $p$ and $a$ be a point in domain of $S$.
Denoted by $\{S\}_a$ the germ of $p$ at $a$. Let $\{S\}_{a_1}$ and
$\{S\}_{a_2}$ be germs of sections $S_1$ and $S_2$ respectively.
We say that $\{S\}_{a_1}$ and $\{S\}_{a_2}$ are $G_+$-equivalent
if there exist $f\in G_+$ such that
$\{f(S_1)\}_{f(a_1)}=\{S\}_{a_2}$. A germ $\{S\}_a$ is
\textit{regular of class i} if there exist a neighborhood $\cal O$
of $a$ and subbundle $E_i$ such that Im$S|_{\cal O}\subset E_i$.
If $\{S\}_a$ is a regular germ of class $i\geq0$, then in a
neighborhood of $a$ we have $S(x)=(x,0,...,0,a_i(x),...,a_0(x))$.
In the rest of the paper we will often denoted $\{S\}_a$ by
$\{a_i,...,a_0\}_a$. If $\{S\}_a$ is a regular germ, then $a$ is a
\textit{regular point} of $S$.
\begin{Definition}
Let $S$ be a section of $p$ and \textbf{v} be a vector field of
the Lie algebra of group $G$, if $\theta_t$ be the flow of
\textbf{v}, we say \textbf{v} is a projective symmetry of S if one
of the following statements satisfied:
\begin{itemize}
\item[1)]$\theta_t(S)=\overline{\theta_t}\circ
S\circ\theta_t^{-1}=S$,
\item[2)]$\displaystyle\frac{d}{dt}\theta_t(S)\Big|_{t=0}=0.$
\end{itemize}
\end{Definition}
Denote by ${\cal P}(S)$ the Lie algebra of all projective
symmetries of $S$.

Let $\Upsilon$ be the set of all regular germs at $0\in\bf R$ of
sections of $p$. Define
\begin{eqnarray*}
\Upsilon_i=\Big\{\{S\}_a|\mbox{dim}{{\cal P}(S)}=i\Big\},\quad
i=0,1,3,
\end{eqnarray*}
and denote $\Upsilon=\Upsilon_0\cup\Upsilon_1\cup\Upsilon_3$. If
$G_0$ be the isotropic subgroup of $G$ in 0, then, $\Upsilon_i$'s
are $G_0$-invariant.

Define $\Upsilon_{r,i}\subset\Upsilon_r$ be the subset of all
regular germs of class $i$. It follows from the invariance of
subbundle $E_i$'s under $G_0$, $\Upsilon_{r,i}$ is
$G_0-$invariant. Consequently we have
\begin{eqnarray*}
\Upsilon_r=\bigcup_{i=0}^{n-3}\Upsilon_{r,i},
\end{eqnarray*}
where this union is separated invariant subsets.

Let ${\bf R}_+$ and ${\bf R}_-$ be the set off positive and
negative real numbers respectively. If
$\ell_{r,i}:\Upsilon_{r,i}\rightarrow({\bf
R}\backslash\{0\})\times{\bf R}$ be a map by the formula
$\{a_i,...,a_0\}\mapsto(a_i(0),a_i'(0))$ and
\begin{eqnarray}\label{eq:11}
&&G_{0+}\times\Upsilon_{r,i}\rightarrow\Upsilon_{r,i}\\
&&(f,\{S\}_0)\mapsto\{f(S)\}_0,\nonumber
\end{eqnarray}
be the action of $G_{0+}$ on $\Upsilon_{r,i}$ then,
\begin{Lemma}\label{lem:2}
The map $\ell_{r,i}|_\Theta$ is a bijection from the orbit
$\Theta$ of the action (\ref{eq:11}) either to $({\bf
R}_+)\times{\bf R}$ or to $({\bf R}_-)\times{\bf R}$.
\end{Lemma}

Let $\Omega_{r,i}^+=\ell_{r,i}^{-1}((1,0))$ and
$\Omega_{r,i}^-=\ell_{r,i}^{-1}((-1,0))$. Denote by $\Gamma_{r,i}$
the subset of $\Omega_{r,i}^+\cup\Omega_{r,i}^-$ defined in the
following way:
\begin{itemize}
\item[1)]$\Gamma_{r,0}=\Omega_{r,0}^+\cup\Omega_{r,0}^-$ for
$i=0$,
\item[2)]if $i>0$, then, $\Omega_{r,i}$ consists of all germs
$\{a_i,...,a_0\}$ from $\Omega_{r,i}^+\cup\Omega_{r,i}^-$
satisfying one of the following conditions:
\begin{itemize}
\item[i)]$a_{i-j}=0$ for all odd numbers $j$ with $1\leq j\leq i$,
\item[ii)]there exist an odd number $r$ with $1\leq r\leq i$ such
that $a_{i-r}(0)>0$ and if $r>1$, then $a_{i-j}(0)=0$ for all odd
numbers $j$ with $1\leq j<r$.
\end{itemize}
\end{itemize}
\subsection{Classification of regular germs from the family $\displaystyle\bf\Omega_{r,i}$}
Let $\mu\in G_-$ defined by $\mu(x)=-x$ for all $x\in\bf R$, then,
due to lemma (\ref{lem:2}) and attentive to
$\mu(\Omega_{r,i}^-)=\Omega_{r,i}^+$ we have:
\begin{Theorem}
\begin{itemize}
\item[1)]The set $\Omega_{r,i}^+\cup\Omega_{r,i}^-$ is a family of
all germs from $\Upsilon_{r,i}$ nonequivalent with respect to
$G_{0+}$.
\item[2)] If $n-i$ is odd, then $\Omega_{r,i}^+$ is a family of
all germs from $\Upsilon_{r,i}$ nonequivalent with respect to
$G_0$.
\item[3)] If $n-i$ is even, $\Gamma_{r,i}$ is a family of
all germs from $\Upsilon_{r,i}$ nonequivalent with respect to
$G_0$.
\end{itemize}
\end{Theorem}
an important corollary concludes this section as follows:
\begin{Corollary}
Classification of regular germs of sections of (\ref{eq:7}) is:
\begin{itemize}
\item[1)] The family of germs of the form
\begin{eqnarray*}
\{\pm1+a(x)x^2,a_{i-1}(x),...,a_0\}_0
\end{eqnarray*}
is a family of all regular germs of class $i$ nonequivalent with
respect to $G_{0+}.$
\item[2)]If $n-i$ is odd, then the family of germs of the form
\begin{eqnarray*}
\{1+a(x)x^2,a_{i-1}(x),...,a_0\}_0
\end{eqnarray*}
is a family of all regular germs of class $i$ nonequivalent with
respect to $G_0.$
\item[3)]If $n-i$ is even, then the family of germs of the form
\begin{eqnarray*}
\{\pm1+a(x)x^2,a_{i-1}(x),...,a_0\}_0,
\end{eqnarray*}
satisfying one of the following conditions:
\begin{itemize}
\item[a)]$a_{i-j}(0)=0$ for all odd numbers $j$ with $1\leq j\leq
i$,
\item[b)]there exist an odd number $r$ with $1\leq r\leq
i$ such that $a_{i-r}(0)>0$ and if $r>1$, then $a_{i-j}(0)=0$ for
all odd number $j$ with $1\leq j\leq r$,
\end{itemize}
is the family of germs of class $i$ nonequivalent with respect to
$G_0.$
\end{itemize}
\end{Corollary}
\section{Conclusion}
This article was a qualification of classification of linear ODEs
due to V.A. Yumaguzhin. First we transform the general form of
ODEs to Laguerre-Forsyth form, then by a suitable change of
variable up to projective transformation we reduce this
classification to classification of the sections of bundles, next
by construction germs and specially regular germs of this sections
near identity, the classification reduced to classifying of
regular germs by providing some invariant subsets of the bundles.

\bibliographystyle{amsplain}
\bibliography{xbib}
\end{document}